\documentclass[a4paper,10pt,DIV=12,flalign, abstracton]{scrartcl}
\usepackage{float}
\usepackage[utf8]{inputenc}
\usepackage[T1]{fontenc}
\usepackage{amssymb,amsmath,amsthm}
\usepackage{graphicx}
\usepackage{hyperref}
\usepackage{csquotes}
\usepackage{subcaption}
\usepackage[sort&compress,numbers]{natbib}
\usepackage{float}
\usepackage{array}

\usepackage[section]{placeins}

\usepackage{color}

\title{An efficient monolithic solution scheme for FE\textsuperscript{2} problems}
\author{Nils Lange, Geralf Hütter, Björn Kiefer}

\begin{document}

\newcolumntype{C}[1]{>{\centering\arraybackslash}m{#1}}

\maketitle
\begin{abstract}
The FE\textsuperscript{2} method is a very flexible but computationally expensive tool for multiscale simulations. In conventional implementations, the microscopic displacements are iteratively solved for within each macroscopic iteration loop, although the macroscopic strains imposed as boundary conditions at the micro-scale only represent estimates. In order to reduce the number of expensive micro-scale iterations, the present contribution presents a monolithic FE\textsuperscript{2} scheme, for which the displacements at the micro-scale and at the macro-scale are solved for in a common Newton-Raphson loop. In this case, the linear system of equations within each iteration is solved by static condensation, so that only very limited modifications to the conventional, staggered scheme are necessary. The proposed monolithic FE\textsuperscript{2} algorithm is implemented into the commercial FE code Abaqus. Benchmark examples demonstrate that the monolithic scheme saves up to $\approx$~60\% of computational costs.


\includegraphics[width=0.7\textwidth]{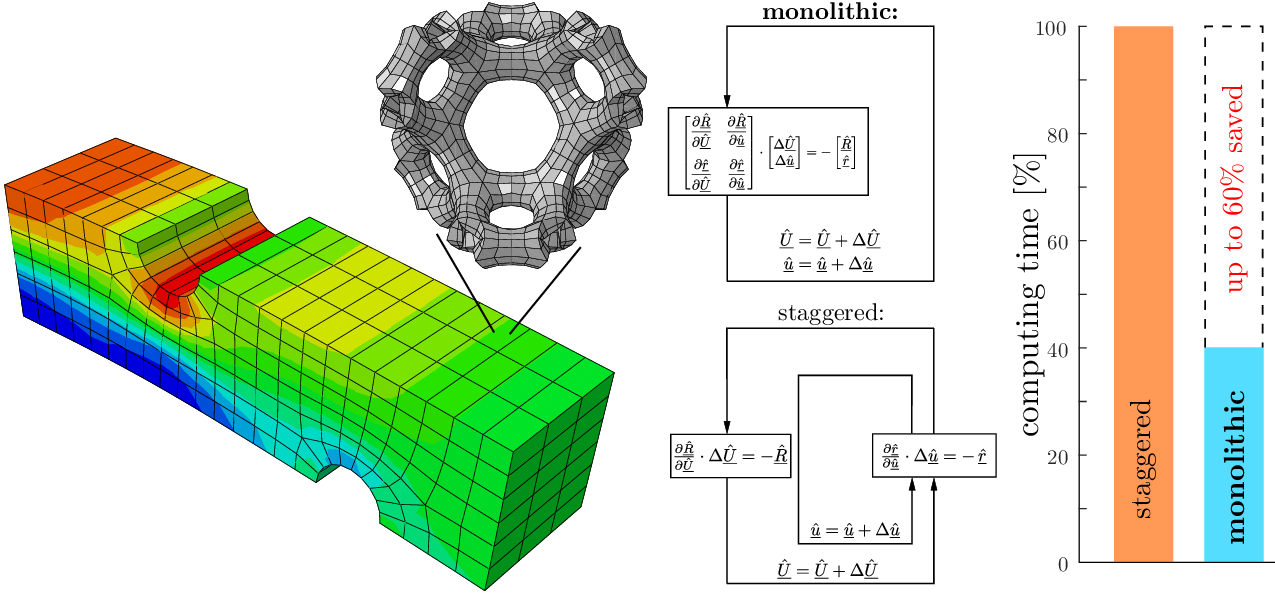}
\end{abstract}

%--------------------------------------------
\section{Introduction}
%--------------------------------------------

\begin{figure}[H]
	\centering
	\includegraphics[width=0.7\linewidth]{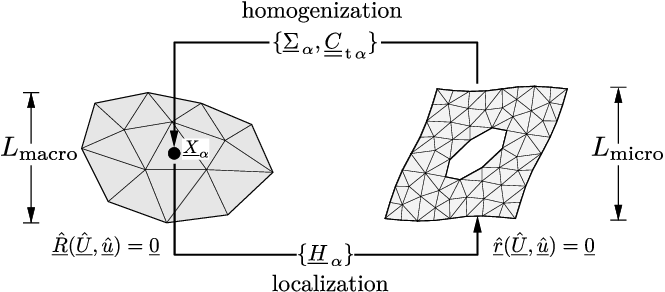}
	\caption[]{Homogenization of the discretized microstructure for a prescribed a macroscopic displacement gradient $\underline{H}_{\,\alpha}$. The solution of the micro-scale boundary value problem associated with the microstructure locally attached at $\underline{X}_{\,\alpha}$ yields the effective macroscopic stress $\underline{\Sigma}_{\,\alpha}$ and tangent stiffness $\underline{\underline{C}}_{\,\text{t}\,\alpha}$ \citep[after]{bayreu_mehrskalen_2005}).}
	\label{fig:FE2_main_principle}
\end{figure}

Most engineering materials are inhomogeneous at a characteristic length-scale. The properties of this so-called micro-scale, consisting of certain dissimilar constituents, determines the behavior at the macro-scale. The macro-scale, in turn, is usually the actual scale of interest (for engineering structures). While it is generally possible to directly resolve the microstructure in structural computations, it is often not favorable, since this results in models far too complex to be dealing with large structures \cite{kanoute_multiscale_2009}. A common way to reduce the model size by means of multi-scale approaches. The mutual connection between the micro- and macro-scales is established through localization and  homogenization processes. In this scale transition, heterogeneous microstructure is assigned to a homogeneous material with effective properties (homogenization), while the inverse procedure is called localization \cite{bayreu_mehrskalen_2005}. A comprehensive overview on the theoretical background of different multi-scale approaches can, for instance, be found in \citep{kanoute_multiscale_2009,Geers_multiscale_2017,Schroeder2014,Raju2021}). 
\\
\noindent The use of finite element analysis (FEA) at both length scales is referred to as FE\textsuperscript{2} modeling \citep{feyel_multiscaleFE2_1999} and illustrated in Figure~\ref{fig:FE2_main_principle}. In linear problems, a sequential solution procedure is possible, since the effective (macro-)stiffness of the microstructure is representative of all deformation states and can therefore be calculated beforehand. For many nonlinear problems, a concurrent solution scheme is necessary, in which a micro-scale model has to be solved at each integration point of the macroscopic mesh \citep{boehm_micromechanics_2018,Schroeder2014,Geers_multiscale_2017}. While the implementation of the sequential scheme is comparatively simple, the concurrent approach is challenging and computationally expensive, but still quite universally applied. To name a few examples, it has been used in numerical multi-scale analysis of fibre-reinforced composites \cite{feyel_multiscaleFE2_1999,Schroeder2014,tikarrouchine_fe2_composite_2018}, woven composites \cite{Yuan_comp_homogenization_2008,tikarrouchine+etal21}, biomechanics \citep{Okada2010}, elastic-plastic matrix-inclusion problems \cite{Miehe2007,Schroeder2014}, or the growth of microvoids \cite{Kouznetsova_micro_macro_2001,Miehe2007} and texture formation in polycrystalline metals \cite{Miehe1999}, and in the design optimization of microstructures \cite{Nakshatrala_multiscale_topology_2013}. Following much research attention over the last two decades, the FE\textsuperscript{2} method is no longer restricted to conventional mechanical problems. Extensions of the method to the realm of multi-physics and generalized continua have been proposed by many groups, e.g.~for thermo-mechanical \cite{Oezdemir_thermo_mech_fe2_2008,berthelsen+menzel19}, electro-mechanical \cite{khalaquzzaman+etal12,keip+steinmann+schroeder14}, magneto-mechanical \cite{zabihyan+etal20}, micromorphic \cite{Biswas2019,Rokos2018} or even three-field problems \cite{schroeder+etal15,labusch+schroeder+keip18} --- see also the references therein. 
\noindent Several implementations can be found in literature. The first FE\textsuperscript{2}-program was programmed by \citet{feyel_multiscaleFE2_1999}. He used the ZeBuLoN code for both scales in the way that the code is reentrant, which means that it is able to call itself. \citet{Kouznetsova_micro_macro_2001} presented an implementation in 2001, where the macroscopic model is computed by a MATLAB FE-code and the microscopic problems are solved by the commercial package MARC. \citet{Yuan_comp_homogenization_2008}, \citet{tchalla_multiscale_abaqus_2013} and \citet{tikarrouchine_fe2_composite_2018,tikarrouchine+etal21} utilize Abaqus for both scales with the aid of a python script and a procedure \texttt{UMAT} to let the Abaqus macro program calls itself at each integration point. 
Naturally this listing is not complete.

The issue of high computational costs has already been addressed in some of the the mentioned publications. \citet{Yuan_comp_homogenization_2008} specify the number of linear solution operations of the microscopic problem with $N_{\text{cells}}\cdot n \cdot I_{\text{macro}}\cdot I_{\text{micro}}$ where  $N_{\text{cells}}$ is the number of macroscopic integration points, $n$ the number of load increments and $I_{\text{macro}}$ and $I_{\text{micro}}$ the average number of iterations on the macro- and micro-scale. As this number can become quite large, even for small problems, the necessity to substantially reduce the computational effort is evident. They name coarse graining or model reduction and parallel computation as the main options to lower computational costs \citep{Yuan_comp_homogenization_2008}. The need for parallel computation is also mentioned in other publications, cf.~\citep{feyel_multiscaleFE2_1999,tikarrouchine_fe2_composite_2018}. Highly parallelizing implementations of FE\textsuperscript{2} have been presented in \citep{Klawonn2020,Mosby2015}.

Some publications suggest algorithmic modifications to the standard FE\textsuperscript{2} solution procedure. Feyel for example mentions the use of a Quasi-Newton algorithm, so that the macroscopic material tangent does not have to be recalculated in every Newton step \cite{feyel_multiscaleFE2_1999}. \citet{temizer+wriggers08} investigate condensation and perturbation procedures in the context of computing the macroscopic tangent. \citet{Nezamabadi_asymptotic_2010} combine the multi-scale finite element procedure with the asymptotic numerical method, which can be efficient for certain types of problems, e.g.~in buckling analysis. {A completely different approach described in the literature is to employ a fast Fourier transformation method to  solve the micro-problem efficiently, while using finite element analysis for the macro-problem, cf.~\citep{Spahn_multiscale_fft_2014,Welschinger_fe_fft_2019, Guodong_fe_fft_2019}.} This spectral approach is, however, limited to voxalized microstructures and its computational costs increase with the contrast between the moduli of the constituents. The idea of solving micromechanical problems by means of Fourier transformations was first developed by Moulinec and Suquet \cite{Moulinec_fft_1994}. An overview over recent applications of this method can be found in \cite{Shanthraj_spectral_Solvers_2019}. Another approach proposed  in several publications towards computationally efficient multi-scale simulations is a sequential solution procedure, which is capable of describing nonlinear and irreversible effects with the aid of neural networks (NN), see \cite{Le_neural_network_2015,Welschinger_fe_fft_2019,settgast_neural_network_2020,abendroth+etal20}. The efficiency of this approach arises from the fact that the constitutive relation at the macroscopic integration point is provided by a trained neural network instead of through concurrent microscopic finite element simulations. Disadvantages are the lower flexibility and the costly training process associated with the NN-approach.

The fundamental equation of the FEM is the following general equilibrium relation, to be fulfilled at both the macro- and micro-levels
\begin{equation}\label{equilibrium_macro}
\underline{\hat{R}}_{ }(\underline{\hat{U}},\ \underline{\hat{u}})=\underline{\hat{F}}_{\,\text{int}}(\underline{\hat{U}},\ \underline{\hat{u}})=\sum_{\alpha=1}^{n_{\alpha}}w_{\alpha}\cdot\underline{\underline{B}}^{T}_{\,\alpha}\cdot\underline{\Sigma}_{\,\alpha}(\underline{\hat{u}}_{\,\alpha},\ \underline{H}_{\,\alpha}=\underline{\underline{B}}_{\,\alpha}\cdot \underline{\hat{U}})=\underline{0}\ ,
\end{equation}
\begin{equation}\label{equilibrium_micro}
\underline{\hat{r}}_{\,\alpha}(\underline{\hat{u}}_{\,\alpha},\ \underline{H}_{\,\alpha}=\underline{\underline{B}}_{\,\alpha}\cdot\underline{\hat{U}})=\underline{0}\ .
\end{equation}
Here small symbols denote microscopic quantities, while large symbols represent a macroscopic quantity. $\underline{\hat{R}},\ \underline{\hat{r}}$ are the residuals, $\underline{\hat{U}},\ \underline{\hat{u}}$ represent nodal displacements and the index $\alpha$ denotes macroscopic integration points. For simplicity, and without the loss of the general applicability of the presented algorithmic concept, only internal forces $\underline{\hat{F}}_{\,\text{int}}$ have been considered in equation \eqref{equilibrium_macro}. It can be seen that \eqref{equilibrium_macro} and \eqref{equilibrium_micro} are coupled through the macroscopic stress $\underline{\Sigma}_{\,\alpha}$ and displacement gradient $\underline{H}_{\,\alpha}$.
\begin{figure}[H]
	\centering
	\includegraphics[width=0.8\linewidth]{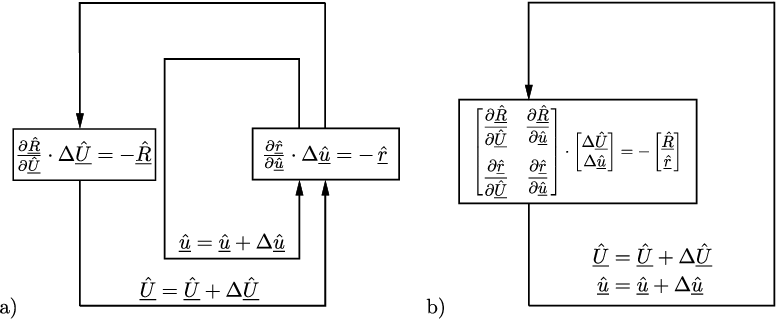}
	\caption{Approaches for solving FE\textsuperscript{2} problems: a) conventional staggered scheme, b) novel monolithic scheme.}
	\label{fig:monolithic_FE2_main_principle}
\end{figure}
The previously cited works all imply a staggered solution scheme, as illustrated in Figure~\ref{fig:monolithic_FE2_main_principle}a). In this context the term \emph{staggered} indicates that a microscopic Newton loop is nested into the macroscopic iterative Newton process. The micro-scale problems receive the macroscopic displacement gradient $\underline{H}_{\,\alpha\,K}$ for the (periodic) boundary condition and iterate until convergence is reached, despite the fact that $\underline{H}_{\,\alpha\,K}$ is only an estimate for $\underline{H}_{\,\alpha}^{\,i+1}$ in the $K$-th Newton-Raphson step. 
Recently, \citet{Tan2020} proposed the \enquote{Direct FE\textsuperscript{2}} method where the macroscopic and microscopic nodes and elements are implemented in a single global FE model, coupled by multi-point constraints. Consequently, the macroscopic and all microscopic problems are solved in a common Newton-Raphson loop in this \emph{monolithic} scheme as shown Figure~\ref{fig:monolithic_FE2_main_principle}b). This method has the advantage that no costly microscopic iterations are necessary and all features of existing FE codes can be used at both scales. However, the global system of equations comprises the macroscopic nodal displacements $\underline{\hat{U}}$ and the microscopic displacements $\underline{\hat{u}}_{\,\alpha}$ of \emph{all} macroscopic Gauss points $\alpha$ and can thus become exceedingly large even for medium-sized problems \citep{Mosby2015}. Interestingly, similar concepts can be found in phase-field fracture modeling approaches, in which efficient monolithic schemes have been implemented and successfully shown to be capable of saving computing time, when compared to the more conventional staggered, algorithmically decoupled approach \cite{borden_phase_field_fracture_2012,gerasimov_monolith_frac_2016}. A more general discussion on this kind of algorithm and generally on iterations at the material level can be found in \citep{Hartmann2005,Yamamoto2020}.

In the present article, a new algorithmic strategy is proposed for solving FE\textsuperscript{2} problems more efficiently in a monolithic way, in which static condensation \citep{Okada2010,temizer+wriggers08} is employed to avoid prohibitively large systems of equations.

%--------------------------------------------
\section{FE\textsuperscript{2}-solution approaches}
%--------------------------------------------

\subsection{Homogenization}
The microscopic FEA can be viewed as a complex materialroutine in the overall workflow of a macroscopic FEM simulation. It replaces the phenomenological material description used in classical FEM. It gives the relation between stress $\underline{\Sigma}_{\,\alpha}$ and displacement gradient $\underline{H}_{\,\alpha}$ \footnote{Kinematic quantities: displacement gradient: $\underline{H}=\underline{U}\otimes\nabla_{\underline{X}_0}$; relation to the deformation gradient: $\underline{F}=\underline{I}+\underline{H}$ and the infinitesimal strain tensor: $\underline{E}=\text{sym}(\underline{H})$.}  in dependence of the load history encountered by the micro model, which is described through internal state variables $\underline{a}$ stored for each microscopic integration point. At the micro-scale a phenomenological material model is employed
\begin{equation}\label{material_law}
\underline{\sigma}=\underline{\sigma}(\underline{h},\underline{a})\ ,
\end{equation}
with appropriate evolution equations for the internal state variables $\underline{a}$.
\\\\
The macroscopic stress $\underline{\Sigma}$ is obtained by homogenizing the microscopic stresses $\underline{\sigma}$ over the volume. It can be shown that the volume integral can be transformed into a surface integral
\begin{equation}\label{homstress}
\underline{\Sigma}=\langle\underline{\sigma}\rangle=\frac{1}{\vert{V_0}\vert}\int_{{V_0}}^{ }\underline{\sigma}\ dV_0=\frac{1}{\vert{V_0}\vert}\int_{\partial{V_0}}^{ }\underline{t}\otimes\underline{x}_0\ dS_0\ .
\end{equation}
To satisfy the condition that the averaged microscopic displacement gradient is equal to the macroscopic displacement gradient
\begin{equation}\label{homE}
\underline{H}=\langle \underline{h}\rangle=\frac{1}{\vert{V_0}\vert}\int_{{V_0}}^{ }\underline{h}\ dV_0=\frac{1}{\vert{V_0}\vert}\int_{\partial{V_0}}^{ }\underline{u}\otimes\underline{n}_0\ dS_0\ ,
\end{equation}
suitable boundary conditions are to be enforced. In this article only periodic boundary conditions are considered, as they estimate the best overall response. They can be formulated as
\begin{equation}\label{PBC}
\underline{u}\left(\underline{x}^+,t\right)=\underline{u}\left(\underline{x}^-,t\right)+\underline{H}(t)\cdot\left[\underline{x}_0^+-\underline{x}_0^-\right]\ .
\end{equation}
The boundary of the microscopic model $\partial{V}$ is split into two parts $\partial{V}=\partial{V}^+\cup\partial{V}^-$. Always two corresponding points on the boundary partitions $\underline{x}^+\in \partial{V}^+$ and $\underline{x}^-\in \partial{V}^-$ with oppositely oriented normal vectors $\underline{n}_0^-=-\,\underline{n}_0^+$ are linked together.
\\\\
Switching to a finite element context, equation \eqref{PBC} has to be enforced through elimination, a penalty parameter, or by using Lagrange multipliers, leading to anti-periodic nodal forces \cite{michel_composite_1999,temizer+wriggers08}. The macroscopic stress is computed in accordance to equation \eqref{homstress} through a summation over the nodal forces of the boundary of the microscopic model. With the macro stress $\underline{\Sigma}_{\,\alpha}$ being a function of the microscopic nodal displacements and the macro displacement gradient
\begin{equation}
\underline{\Sigma}_{\,\alpha}=\underline{\Sigma}(\underline{\hat{u}}_{\,\alpha},\underline{H}_{\,\alpha})\ ,
\end{equation}
the macroscopic material tangent $\underline{\underline{C}}_{\,\text{t}\,\alpha}$ can be calculated as
\begin{equation}\label{mat_tangent1}
\underline{\underline{C}}_{\,\text{t}\,\alpha}=\frac{\text{d}\underline{\Sigma}_{\alpha}}{\text{d}\underline{H}_{\,\alpha}}=\frac{\partial\underline{\Sigma}_{\,\alpha}}{\partial\underline{H}_{\,\alpha}}+\frac{\partial\underline{\Sigma}_{\,\alpha}}{\partial\underline{\hat{u}}_{\,\alpha}}\cdot\frac{\partial\underline{\hat{u}}_{\,\alpha}}{\partial\underline{H}_{\,\alpha}}\ .
\end{equation}
Here, $\frac{\partial\underline{\hat{u}}_{\,\alpha}}{\partial\underline{H}_{\,\alpha}}$ is only implicitly known over equation \eqref{equilibrium_micro}. An implicit differentiation yields
\begin{equation}\label{mat_tangent2}
\underbrace{\dfrac{\partial \underline{\hat{r}}_{\,\alpha}}{\partial \underline{\hat{u}}_{\,\alpha}}}_{=:\underline{\underline{k}}_{\,\text{t}\,\alpha}}\cdot\frac{\partial\underline{\hat{u}}_{\,\alpha}}{\partial\underline{H}_{\,\alpha}}+\frac{\partial\underline{\hat{r}}_{\,\alpha}}{\partial\underline{H}_{\,\alpha}}=\underline{0}\ \leftrightarrow\ \frac{\partial\underline{\hat{u}}_{\,\alpha}}{\partial\underline{H}_{\,\alpha}}=-\left[\underline{\underline{k}}_{\,\text{t}\,\alpha}\right]^{-1}\cdot\frac{\partial\underline{\hat{r}}_{\,\alpha}}{\partial\underline{H}_{\,\alpha}}\ .
\end{equation}
Therein, $\underline{\underline{k}}_{\,\text{t}\,\alpha}$ denotes the tangent stiffness matrix of the microscopic problem connected to a macroscopic integration point $\alpha$. Inserting equation \eqref{mat_tangent2} into \eqref{mat_tangent1}, the macroscopic tangent stiffness can formally be written as
\begin{equation}\label{mat_tangent3}
\underline{\underline{C}}_{\,\text{t}\,\alpha}=\frac{\text{d}\underline{\Sigma}_{\,\alpha}}{\text{d}\underline{H}_{\,\alpha}}=\frac{\partial\underline{\Sigma}_{\,\alpha}}{\partial\underline{H}_{\,\alpha}}-\frac{\partial\underline{\Sigma}_{\,\alpha}}{\partial\underline{\hat{u}}_{\,\alpha}}\cdot\left[\underline{\underline{k}}_{\,\text{t}\,\alpha}\right]^{-1}\cdot\frac{\partial\underline{\hat{r}}_{\,\alpha}}{\partial\underline{H}_{\,\alpha}}\ .
\end{equation}
Though, the inverse in equation \eqref{mat_tangent3} is not calculated in the actual implementation. Instead, the system of equations \eqref{mat_tangent2} is solved for $\frac{\partial\underline{\hat{u}}_{\,\alpha}}{\partial\underline{H}_{\,\alpha}}$.\\\\
When small deformations are considered, it is sufficient to prescribe the infinitesimal macroscopic strain tensor $\underline{E}=\text{sym}(\underline{H})$ in equation \eqref{PBC} and the stress $\underline{\Sigma}$ corresponds to the Cauchy stress. In the finite deformation theory $\underline{\Sigma}$ corresponds to the 1. Piola Kirchhoff stress \cite{hill_1984}.
\subsection{Staggered solution scheme}
The discretized equilibrium conditions \eqref{equilibrium_macro} and \eqref{equilibrium_micro} represent a nonlinear system of equations. They have to be solved utilizing a numerical method, generally the Newton Raphson (NR) method or its modifications. In the conventional, staggered solution procedure the macroscopic and all microscopic problems are treated separately. It results in nesting the microscopic NR loop into the macroscopic one as shown in Figure \ref{fig:monolithic_FE2_main_principle} a). Linearization of equations \eqref{equilibrium_macro} and \eqref{equilibrium_micro} yields
\begin{equation}\label{linearized_micro_equilibrium_staggered}
\underline{\underline{k}}_{\,\text{t}\,\alpha\,k}^{\,K}\cdot\Delta\underline{\hat{u}}_{\,\alpha\,k}^{\,K}=-\,\underline{\hat{r}}_{\,\alpha\,k}^{\,K}\ ,
\end{equation}
\begin{equation}
\underline{\hat{u}}_{\,\alpha\,k+1}^{\,K}=\underline{\hat{u}}_{\,\alpha\,k}^{\,K}+\Delta\underline{\hat{u}}_{\,\alpha\,k}^{\,K}
\end{equation}
and
\begin{equation}
\underline{\underline{K}}_{\,\text{t}\,K}\cdot\Delta\underline{\hat{U}}_{\,K}=-\,\underline{\hat{R}}_{\,K}\ ,
\end{equation}
\begin{equation}
\underline{\hat{U}}_{\,K+1}=\underline{\hat{U}}_{\,K}+\Delta\underline{\hat{U}}_{\,K}\ .
\end{equation}
Here $K$ and $k$ denotes the macroscopic resp.\ microscopic NR index. In the $K$-th macroscopic NR step, the microscopic problems are iterated until convergence is reached. The macroscopic NR procedure continues when the microscopic problems of all macroscopic integration points are equilibrated.
\subsection{Monolithic solution scheme}\label{sec:monolithic_solution_scheme}
In this paper a new algorithmic concept is proposed. Since equations \eqref{equilibrium_macro} and \eqref{equilibrium_micro} are coupled through the macroscopic stress $\underline{\Sigma}_{\,\alpha}$ and displacement gradient $\underline{H}_{\,\alpha}$, they are linearized together
\begin{equation}\label{monolithic_linearization}
\begin{bmatrix}
\dfrac{\partial \underline{\hat{R}}}{\partial \underline{\hat{U}}} & \dfrac{\partial \underline{\hat{R}}}{\partial \underline{\hat{u}}}\\[.5cm]
\dfrac{\partial \underline{\hat{r}}}{\partial \underline{\hat{U}}} & \dfrac{\partial \underline{\hat{r}}}{\partial \underline{\hat{u}}}
\end{bmatrix}
\cdot \begin{bmatrix}
\Delta\underline{\hat{U}}\\ 
\Delta\underline{\hat{u}}
\end{bmatrix}
= - \begin{bmatrix}
\underline{\hat{R}}\\ 
\underline{\hat{r}}
\end{bmatrix}
\end{equation}
The constituents of the overall stiffness matrix of \eqref{monolithic_linearization} can be identified as
\begin{equation}\label{d_R_d_U}
\frac{\partial \underline{\hat{R}}}{\partial \underline{\hat{U}}}=
\sum_{\alpha=1}^{n_{\alpha}} w_{\alpha}\cdot\underline{\underline{B}}^{T}_{\,\alpha}\cdot\underline{\underline{C}}_{\,\text{t}\,\alpha}\cdot\underline{\underline{B}}_{\,\alpha}\ ,
\end{equation}
\begin{equation}\label{d_R_d_u}
\frac{\partial \underline{\hat{R}}}{\partial \underline{\hat{u}}}=
\sum_{\alpha=1}^{n_{\alpha}}w_{\alpha}\cdot\underline{\underline{B}}^{T}_{\,\alpha}\cdot\frac{\partial \underline{\Sigma}_{\,\alpha}}{\partial \underline{\hat{u}}_{\,\alpha}}\,
\end{equation}
and
\begin{equation}\label{d_r_d_U}
\frac{\partial \underline{\hat{r}}}{\partial \underline{\hat{U}}}=\sum_{\alpha=1}^{n_{\alpha}} \frac{\partial \underline{\hat{r}}_{\,\alpha}}{\partial \underline{H}_{\,\alpha}}\cdot\underline{\underline{B}}_{\,\alpha}\ .
\end{equation}
The linear system is generally huge as it comprises the microscopic displacement vector $\underline{\hat{u}}$ for \emph{all} macroscopic integration points. Thus, it is often hardly possible being solved within a single step. As the microscopic nodal displacements are only weakly coupled with the overall system over the macroscopic displacement gradient $\underline{H}_{\,\alpha}$, they can be eliminated through static condensation. Beginning with the second line of \eqref{monolithic_linearization}
\begin{equation}
\frac{\partial \underline{\hat{r}}}{\partial \underline{\hat{U}}}\cdot\Delta\underline{\hat{U}}+\frac{\partial \underline{\hat{r}}}{\partial \underline{\hat{u}}}\cdot\Delta\underline{\hat{u}}= -\,\underline{\hat{r}}\ ,
\end{equation}
inserting \eqref{d_r_d_U} and breaking it down to a summation over the macroscopic integration points yields
\begin{equation}
\sum_{\alpha=1}^{n_{\alpha}}\frac{\partial \underline{\hat{r}}_{\,\alpha}}{\partial \underline{H}_{\,\alpha}} \cdot\underbrace{\underline{\underline{B}}_{\,\alpha}\cdot\Delta\underline{\hat{U}}}_{\Delta\underline{H}_{\,\alpha}}+\sum_{\alpha=1}^{n_{\alpha}}\underbrace{\frac{\partial \underline{\hat{r}}_{\,\alpha}}{\partial \underline{\hat{u}}_{\,\alpha}}}_{\underline{\underline{k}}_{\,\text{t}\,\alpha}}\cdot\Delta\underline{\hat{u}}=- \sum_{\alpha=1}^{n_{\alpha}}\underline{\hat{r}}_{\,\alpha}\ .
\end{equation}
It shows that the (microscopic) systems of equations of all macroscopic integration points are fully decoupled from another
\begin{equation}\label{linearization_monolith_line2}
\sum_{\alpha=1}^{n_{\alpha}}\underbrace{\left[\frac{\partial \underline{\hat{r}}_{\,\alpha}}{\partial \underline{H}_{\,\alpha}}\cdot\Delta\underline{H}_{\,\alpha}+\underline{\underline{k}}_{\,\text{t}\,\alpha}\cdot\Delta\underline{\hat{u}}_{\,\alpha}+\underline{\hat{r}}_{\,\alpha}\right]}_{\underline{0}}=\underline{0}\ .
\end{equation}
Solving this equation for the increment of the microscopic nodal displacement vector yields
\begin{equation}\label{displacement_update_micro_monolith}
\Delta \underline{\hat{u}}_{\,\alpha}=-\left[\underline{\underline{k}}_{\,\text{t}\,\alpha}\right]^{-1}\cdot\left[\underline{\hat{r}}_{\,\alpha}+\frac{\partial \underline{\hat{r}}_{\,\alpha}}{\partial \underline{H}_{\,\alpha}}\cdot\Delta\underline{H}_{\,\alpha}\right]\ .
\end{equation}
Finally the first line of \eqref{monolithic_linearization} is considered
\begin{equation}
\frac{\partial \underline{\hat{R}}}{\partial \underline{\hat{U}}}\cdot\Delta\underline{\hat{U}}+\frac{\partial \underline{\hat{R}}}{\partial \underline{\hat{u}}}\cdot\Delta\underline{\hat{u}}= -\,\underline{\hat{R}}
\end{equation}
and the linearized macroscopic equilibrium equation is obtained with use of \eqref{d_R_d_U}, \eqref{d_R_d_u} and \eqref{equilibrium_macro} as
\begin{equation}\label{monolith_macro_equilib}
\sum_{\alpha=1}^{n_{\alpha}} w_{\alpha}\cdot\underline{\underline{B}}^{T}_{\,\alpha}\cdot\underbrace{\left[\frac{\partial \underline{\Sigma}_{\,\alpha}}{\partial \underline{H}_{\,\alpha}}-\frac{\partial \underline{\Sigma}_{\,\alpha}}{\partial \underline{u}_{\,\alpha}}\cdot\left[\underline{\underline{k}}_{\,\text{t}\,\alpha}\right]^{-1}\cdot\frac{\partial \underline{\hat{r}}_{\,\alpha}}{\partial \underline{H}_{\,\alpha}}\right]}_{:=\underline{\underline{C}}_{\,\text{t}\,\alpha}}\cdot\underline{\underline{B}}_{\,\alpha}\cdot\Delta\underline{\hat{U}}=-\sum_{\alpha=1}^{n_{\alpha}} w_{\alpha}\cdot\underline{\underline{B}}^{T}_{\,\alpha}\cdot\underbrace{\left[\underline{\Sigma}_{\,\alpha}-\frac{\partial \underline{\Sigma}_{\,\alpha}}{\partial \underline{u}_{\,\alpha}}\cdot\left[\underline{\underline{k}}_{\,\text{t}\,\alpha}\right]^{-1}\cdot\underline{\hat{r}}_{\,\alpha}\right]}_{:=\underline{\Sigma}^{\,\text{alg}}_{\,\alpha}}\ .
\end{equation}
Identifying the tangent stiffness
\begin{equation}\label{macro_tangent_stiffness}
\underline{\underline{C}}_{\,\text{t}\,\alpha}=\frac{\partial \underline{\Sigma}_{\,\alpha}}{\partial \underline{H}_{\,\alpha}}-\frac{\partial \underline{\Sigma}_{\,\alpha}}{\partial \underline{\hat{u}}_{\,\alpha}}\cdot\left[\underline{\underline{k}}_{\,\text{t}\,\alpha}\right]^{-1}\cdot\frac{\partial \underline{\hat{r}}_{\,\alpha}}{\partial \underline{H}_{\,\alpha}}
\end{equation}
and the (algorithmic consistent) stress in \eqref{monolith_macro_equilib} as
\begin{equation}\label{macro_alogrithm_stress}
\underline{\Sigma}^{\,\text{alg}}_{\,\alpha}=\underline{\Sigma}_{\,\alpha}-\frac{\partial \underline{\Sigma}_{\,\alpha}}{\partial \underline{\hat{u}}_{\,\alpha}}\cdot\left[\underline{\underline{k}}_{\,\text{t}\,\alpha}\right]^{-1}\cdot\underline{\hat{r}}_{\,\alpha}
\end{equation}
shows that the (algorithmically consistent) tangent stiffness of the monolithic approach does not change in comparison to the staggered algorithm. However, the algorithmic consistent stress $\underline{\Sigma}_{\,\alpha}^{\,\text{alg}}$ comprises a contribution of the residual $\underline{\hat{r}}_{\,\alpha}$ at the micro-scale in addition to the actual (estimate of the) macroscopic stress $\underline{\Sigma}_{\,\alpha}$.

\begin{figure}[!b]
	\centering
	\includegraphics[width=\linewidth]{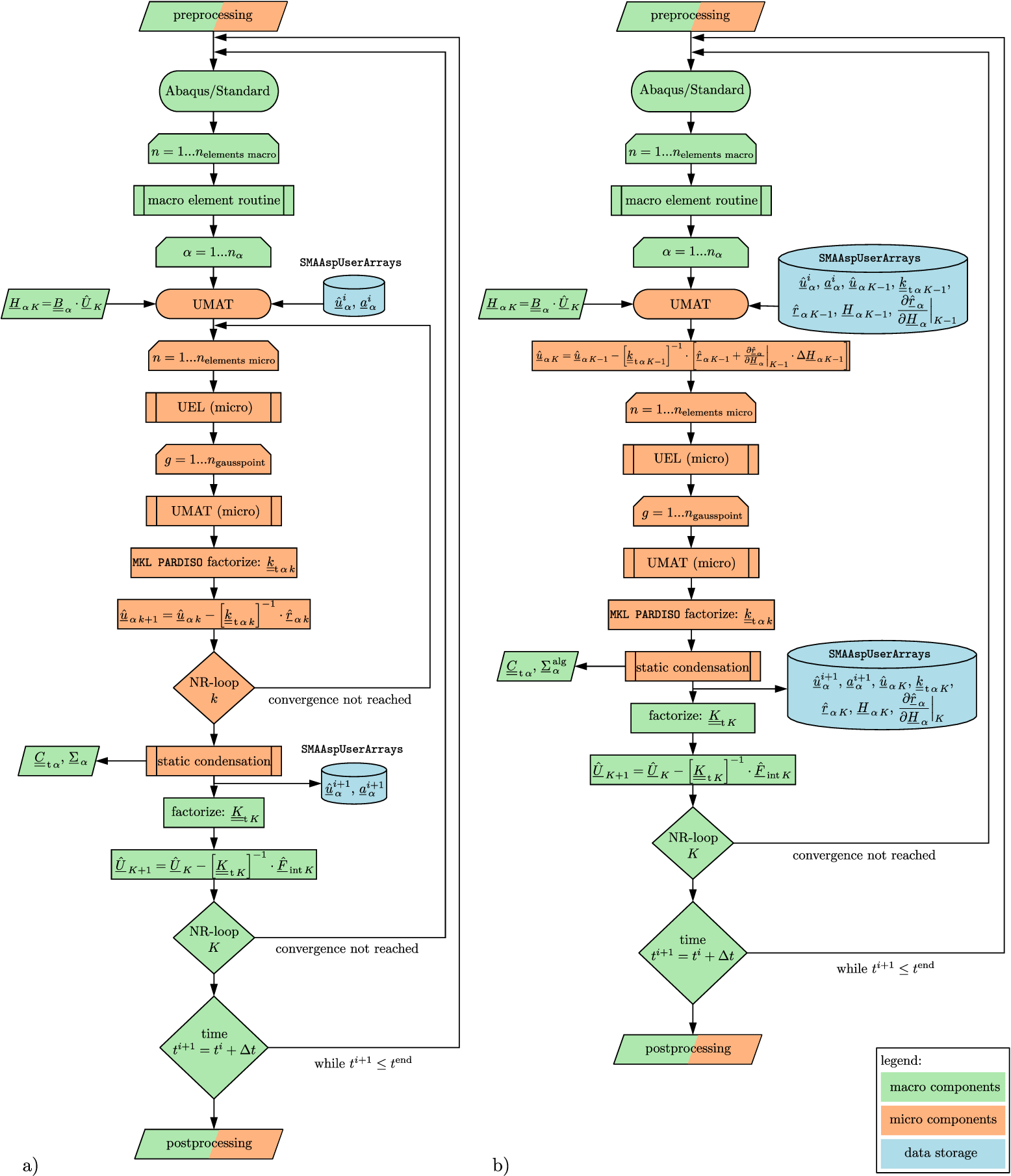}
	\caption{FE\textsuperscript{2} program flowcharts: a) staggered FE\textsuperscript{2} program flowchart; b) monolithic FE\textsuperscript{2} program flowchart.}
	\label{fig:program_flow_charts}
\end{figure}

\subsection{Comparison}
Section \ref{sec:monolithic_solution_scheme} showed, that the monolithic approach can be seamlessly integrated in existing FE code with the following adaptions:
\begin{itemize}
	\item no iteration on the micro level
	\item additional term in the macroscopic stress $\underline{\Sigma}_{\,\alpha}^{\,\text{alg}}$ according to \eqref{macro_alogrithm_stress}
	\item additional term in calculating the increment of the microscopic nodal displacements according to \eqref{displacement_update_micro_monolith}
	\item microscopic FE program must have access not only to the values of the last converged NR step but also of the ones of the current NR step
\end{itemize}
Figure \ref{fig:program_flow_charts} shows the program flow chart of the staggered algorithm on the left and the monolithic on the right. It is to be noted that the evaluation of equation \eqref{displacement_update_micro_monolith} in the monolithic algorithm at the beginning of the next NR increment leads to a straightforward implementation and is algorithmically equivalent to an evaluation at the end of the current NR increment.

%--------------------------------------------
\section{Implementation}
%--------------------------------------------

The monolithic as well as the staggered FE\textsuperscript{2} method have been implemented into the commercial FE code Abaqus.
Existing Abaqus-FE\textsuperscript{2} implementations are build on the idea to use Abaqus both on the macro- and microscopic level, see~\cite{Yuan_comp_homogenization_2008,tchalla_multiscale_abaqus_2013,tikarrouchine_fe2_composite_2018}. The advantage is, that Abaqus is already well tested and has a lot of built-in features, which can be used for the micro and macro model. However it has the big disadvantage, that the Abaqus code is not reentrant and has to be restarted via a python script.
This restart firstly requires time-consuming reading and writing data from and to the hard-drive. Secondly, the mandatory license check in each of the restarted simulations requires additional time.

In order to avoid these time consuming operations, the main idea behind the present implementation is to use Abaqus at the macroscopic scale together with a light-weight self-written (\texttt{FORTRAN 90}) code on the microscopic scale, 
both, for the monolithic as well as the staggered FE\textsuperscript{2} algorithms.
The data exchange between the microscopic and macroscopic scales is performed through the \texttt{UMAT} interface of Abaqus for user-defined materials \cite{abaqus}.
It provides the (macroscopic) strain or deformation gradient, respectively, at the element integration points and expects  the Cauchy stress and the corresponding, consistent tangent stiffness as output.

As the provided Abaqus data interface for internal state variables \texttt{STATEV} is limited to $10,000$ float numbers and only provides data from the last converged time step, a self-written data management had to be implemented. This data management employs thread-safe allocatable arrays through the \texttt{SMAAspUserArrays} interface of Abaqus \cite{abaqus}.

The present program itself calls its element and material routines for the microscale through interfaces equivalent to Abaqus' UEL und UMAT, respectively. This makes the element and material routines in the code compatible with Abaqus.

In order to obtain a faster and more robust convergence of the Newton-Raphson scheme, the initial guesses of the nodal displacements are  linearly extrapolated from the two previous time increments, both, at the microscopic and the macroscopic level.

The linear systems of equations at the microscopic scale, equations \eqref{linearized_micro_equilibrium_staggered}, \eqref{displacement_update_micro_monolith}, \eqref{macro_tangent_stiffness}, \eqref{macro_alogrithm_stress}, are solved by the Intel MKL PARDISO solver, a direct sparse matrix solver \cite{intel_mkl}. The matrices are stored in the Compressed Sparse Row (CSR) Format. 
A matrix is factorized before it can be solved for different right-hand sides.
The monolithic program can make use of that feature, since the factorization $\underline{\underline{k}}_{\,\text{t}\,\alpha}$ that is needed at the end of the program to calculate $\underline{\Sigma}^{\,\text{alg}}_{\,\alpha}$ and $\underline{\underline{C}}_{\,\text{t}\,\alpha}$ is needed again at the beginning of the next NR iteration to compute $\Delta \underline{\hat{u}}_{\,\alpha}$ according to \eqref{displacement_update_micro_monolith}. As the factorization process is the computationally most costly step, this approach can save considerable amount of computing time, though, at the cost of higher memory consumption. In the following this is referred to as ``monolithic scheme w/ stored factorization''.
In an alternative ``monolithic implementation without stored factorization'', the factorization of $\underline{\underline{k}}_{\,\text{t}\,\alpha}$ is not saved but computed again for the update of  $\Delta \underline{\hat{u}}_{\,\alpha}$, resulting in lower memory consumption, but higher computational cost.

In order to reduce the bandwidth $b$ of $\underline{\underline{k}}_{\,\text{t}\,\alpha}$, a permutation matrix is once calculated at the beginning of the simulation, which reorders the degrees of freedom.

The implementation is fully parallelizable in the sense that microscopic problems can be run in parallel. In addition, the macroscopic as well as the microscopic solver can work in parallel mode.\footnote{The program is accessible at \url{https://tu-freiberg.de/MonolithFE2}.}

%--------------------------------------------
\section{Benchmarking}
%--------------------------------------------

To evaluate differences in the computational effort, benchmarks between the staggered and the monolithic (w \& w/o stored factorized stiffness matrix) algorithm are performed, by computing model problems with highly nonlinear character. Thereby the computing time strongly depends on the nonlinearity of the problem chosen. In linear problems the staggered and monolithic algorithm are equivalent. To get comparable results, all examples were calculated on the same computer with one processor (\texttt{Intel(R) Xeon(R) Gold 6244 CPU @ 3.60GHz}). An automatic time increment control was used in Abaqus. When one or more microscopic iterations did not converge in a given number $n_{\text{max}}$ of increments in the staggered mode, Abaqus is directed to lower the time increment, vice versa when all microscopic iterations converged in less or equal $n_{\text{max}}\cdot0.5$ iterations, Abaqus is directed to increase the step size.
\begin{figure}[!h]
	\centering
	\includegraphics[width=0.8\linewidth]{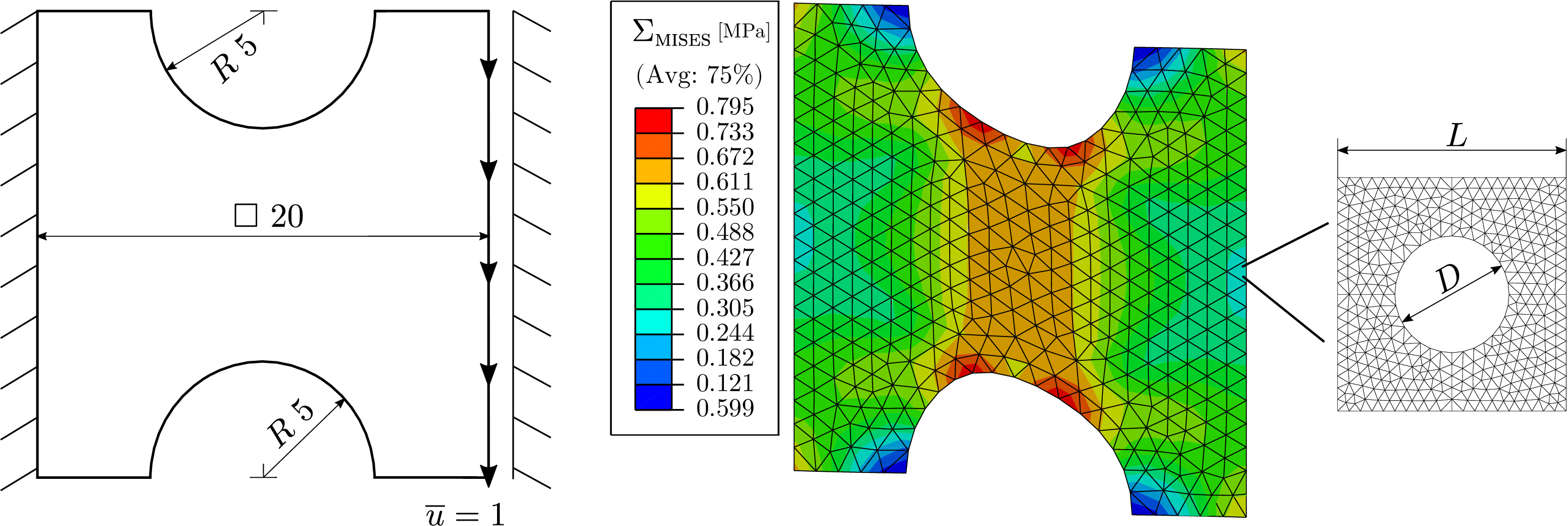}
	\caption{Example 1: notched plate under shear, porous microstructure: $E=100.0$ MPa, $\nu=0.3$, $\sigma_0=1.0$ MPa, $h=2.0$ MPa, plain strain state, $D/L=0.5$, small deformation theory.}
	\label{fig:Benchmarking_notched_plate_shear_problem}
\end{figure}

\begin{figure}[!htb]
	\centering
	\begin{subfigure}[c]{0.49\linewidth}
		\centering
		\captionsetup{justification=centering,margin=0.1\linewidth}
		\includegraphics[width=0.7\linewidth]{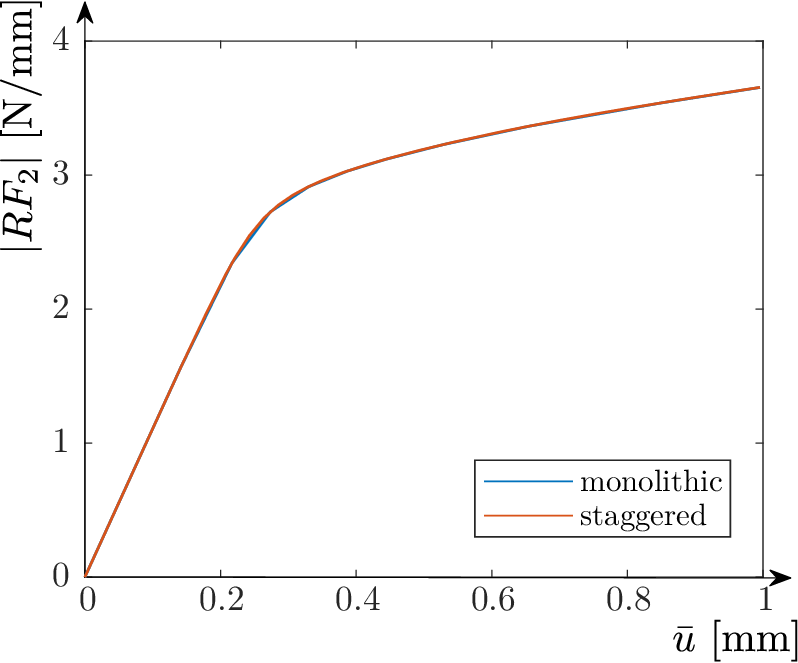}
		\subcaption{total reaction force over prescribed displacement benchmark example 1}
		\label{fig:U_RF_notched_plate_shear}
	\end{subfigure}
	\begin{subfigure}[c]{0.49\linewidth}
		\centering
		\captionsetup{justification=centering,margin=0.1\linewidth}
		\includegraphics[width=0.7\linewidth]{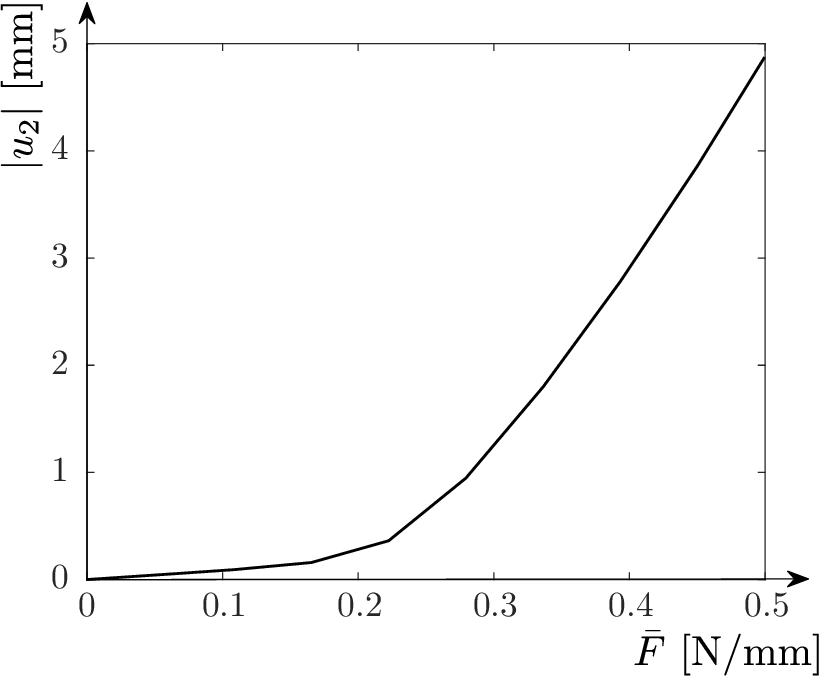}
		\subcaption{horizontal displacement of load application node over prescribed force benchmark example 2}
		\label{fig:F_U_notched_beam_bending}
	\end{subfigure}
	\begin{subfigure}[c]{0.49\linewidth}
		\centering
		\captionsetup{justification=centering,margin=0.1\linewidth}
		\includegraphics[width=0.7\linewidth]{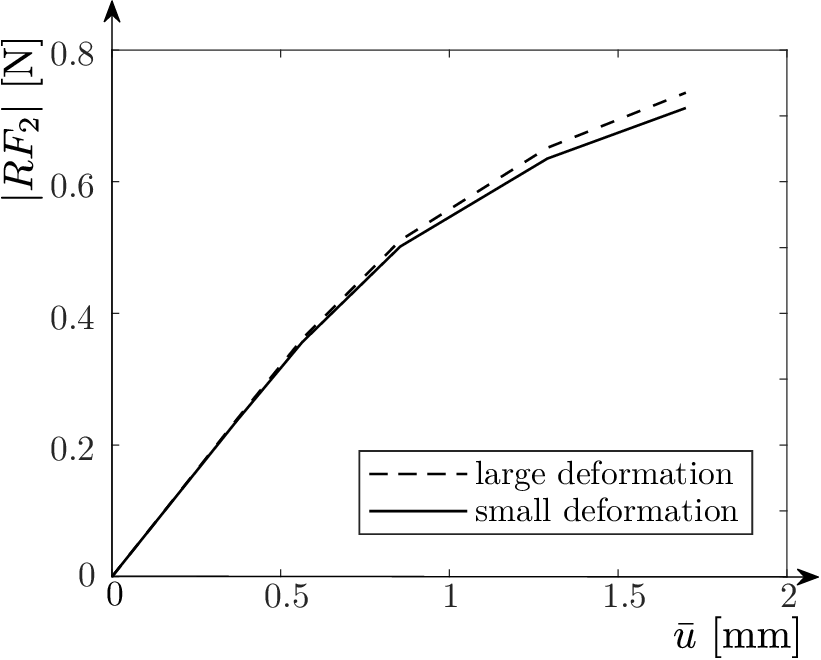}
		\subcaption{total reaction force over prescribed displacement benchmark example 3}
		\label{fig:U_RF_notched_3D_beam_bending}
	\end{subfigure} 
	\begin{subfigure}[c]{0.49\linewidth}
		\centering
		\captionsetup{justification=centering,margin=0.1\linewidth}
		\includegraphics[width=0.7\linewidth]{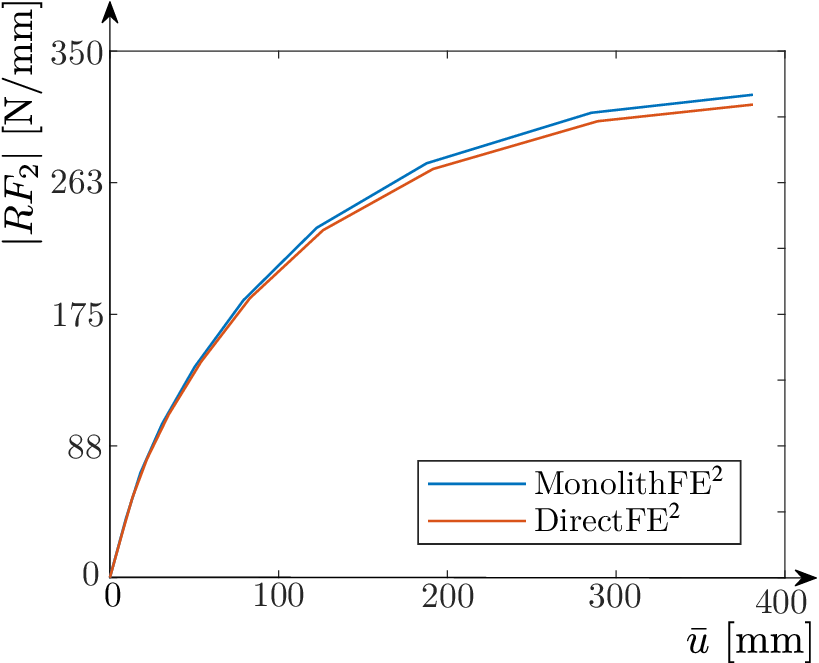}
		\subcaption{total reaction force over prescribed displacement example 4}
		\label{fig:U_RF_beam_direct_fe2}
	\end{subfigure}
	\caption{Global load-displacement curves of all benchmark examples.}
\end{figure}

In the first example a notched plate is exposed to a shear like loading. The macroscopic model, as shown in Figure \ref{fig:Benchmarking_notched_plate_shear_problem}, is meshed with 730 triangular, linear elements. The microscopic model is a porous material having a elastic-M\textsc{ises}-plastic material with linear hardening (initial yield stress $\sigma_0$, hardening modulus $h$). The mesh consists of 706 triangular, linear elements. The total number of degrees of freedom is $538,041$. Figure \ref{fig:U_RF_notched_plate_shear} illustrates the nonlinearity, by plotting the total reaction force over the prescribed displacement for the staggered and monolithic algorithm. It can be seen that both algorithms yield numerically virtually identical global reaction forces. The minor differences arise from the different iteration patterns. The monolithic algorithm saves about 38\% of the computational costs compared to the staggered scheme.
\begin{figure}[!h]
	\centering
	\includegraphics[width=\linewidth]{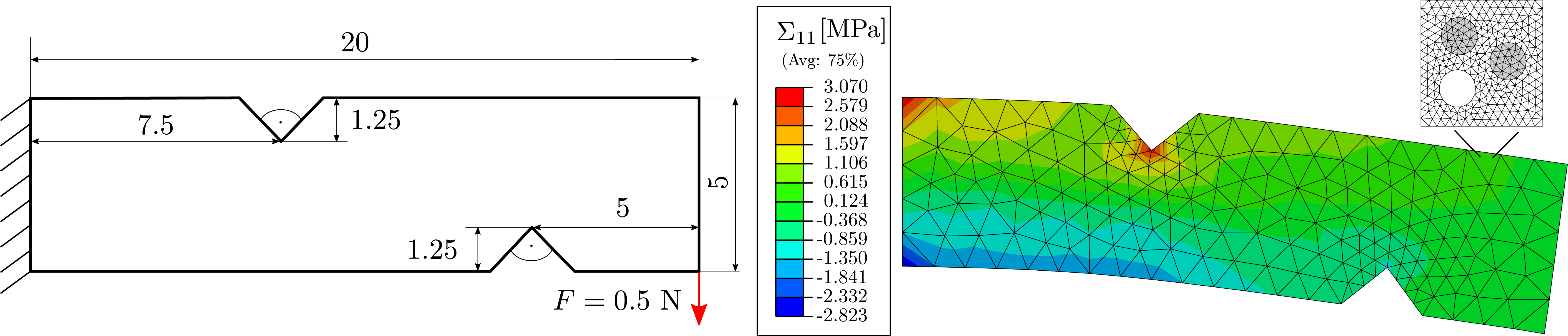}
	\caption{Example 2: notched (2D-)beam undergoing bending, composite microstructure taken from \cite{miehekoch_mikro_makro_2002}, small deformation theory.}
	\label{fig:Benchmarking_notched_beam_bending_problem}
\end{figure}

As second example a two-dimensional, notched beam is exposed to bending. Figure \ref{fig:Benchmarking_notched_beam_bending_problem} shows the macroscopic model, which is meshed with 337 triangular, linear elements. The microscopic model is a composite material with a elastic-plastic matrix, a pore and two stiff, elastic inclusions, taken from \cite{miehekoch_mikro_makro_2002}. The mesh consists of 559 triangular, linear elements. The total number of degrees of freedom is $193,826$. To demonstrate the nonlinearity, Figure \ref{fig:F_U_notched_beam_bending} plots the displacement of the load application node over the prescribed force. The monolithic algorithm reduces the computational effort by 60\% compared to the staggered scheme.
\begin{figure}[!h]
	\centering
	\includegraphics[width=\linewidth]{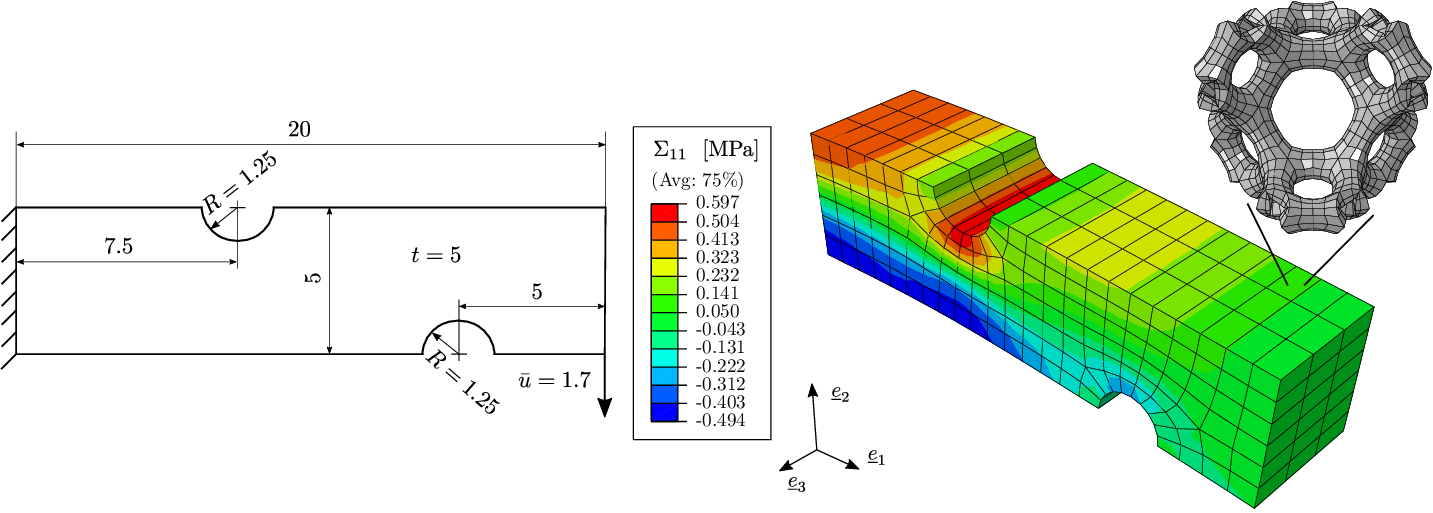}
	\caption{Example 3: notched (3D-)beam undergoing bending, foam microstructure with elastic-(M\textsc{ises})-plastic behavior (hardening after Ramberg–Osgood relation, $E=1630.0$ MPa,\ $\sigma_0=10.0$ MPa,\ $n=10$,\ $\alpha=0.326$), large deformation theory.}
	\label{fig:Benchmarking_notched_3D_beam_bending_problem}
\end{figure}

In the third example a three-dimensional, notched cantilever beam is meshed with 462 linear, hexahedral elements, as shown in Figure \ref{fig:Benchmarking_notched_3D_beam_bending_problem}. The microscopic model from \citep{abendroth_ceramic_foam_2017} is a foam with elastic-plastic behavior and power law hardening according to the Ramberg–Osgood relationship:
\begin{equation}
\varepsilon=\dfrac{\sigma}{E}+\alpha\cdot\dfrac{\sigma}{E}\cdot\left[\dfrac{\sigma}{\sigma_0}\right]^{n-1}
\end{equation}
The mesh of the RVE consists of 2304 linear, hexahedral elements. The foam FE model is descended from the collaborative research center SFB 920 \cite{abendroth_ceramic_foam_2017}. The total number of degrees of freedom amounts to $4,691,042$. The monolithic algorithm saves around 50\% of the computational effort compared to the staggered scheme. The nonlinear character is shown in Figure \ref{fig:U_RF_notched_3D_beam_bending}, where additionally the result for a respective large deformation analysis is plotted.

As it is interesting to investigate how the different program modules (Figure \ref{fig:program_flow_charts}) contribute to the total computing time of the RVE calculations, these times were investigated for the third example. The four main steps are calling all element routines and assembling the stiffness matrices, factorizing the stiffness matrices, solving the system for the increment of the nodal displacements and computing the macroscopic stress and stiffness by static condensation. The computing time of these phases show a ratio of approximately 10:10:1:7 for the third example. Of course, the ratio might change, e.g., if more expensive material routines are used.

\begin{figure}[!h]
	\centering
	\includegraphics[width=0.85\linewidth]{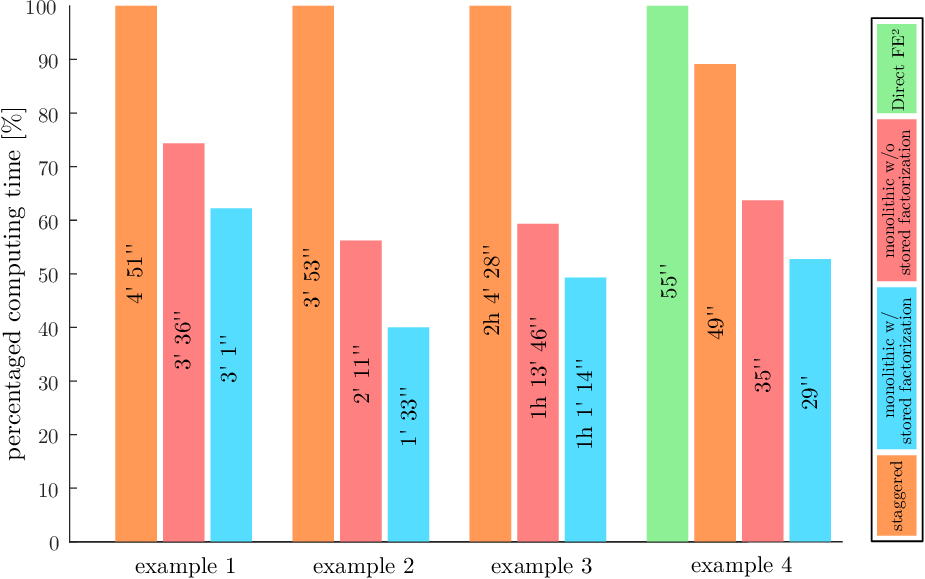}
	\caption{Comparison of benchmarking results (small deformation theory).}
	\label{fig:all_benchmarking_results}
\end{figure}

To get an impression of the memory consumption, Table~\ref{memory} lists the permanent memory needed per macroscopic integration point. The total amount of memory corresponds to the number of integration points times the memory specified. Note that in comparison the memory consumption of the macroscopic model can be neglected. It can be seen that the monolithic implementation without stored factorization requires approximately 1.5 to 2 times as much memory as the staggered scheme. If the factorization is stored, the memory comsumption increases by an additional factor of about 5. Thus, the aspect of memory consumption has to be kept in mind when choosing an algorithm for solving a particular problem on a specific machine. Relating the memory consumption in Table~\ref{memory} with the attainable speed-up factors in Figure~\ref{fig:all_benchmarking_results} shows that the monolithic scheme without storage of the factorized matrices is a reasonable trade-off if memory limitations become relevant.

\begin{table}
	\centering
	\begin{tabular}{C{0.12\linewidth}C{0.22\linewidth}C{0.22\linewidth}C{0.22\linewidth}}
		\hline 
		\textbf{algorithm}& \textit{staggered} & \textit{monolithic w/o stored factorization} & \textit{monolithic w/ stored factorization} \\ 
		\hline 
		\textbf{example 1}& 132.6 & 194.3 & 1078.0 \\ 
		\hline 
		\textbf{example 2}& 103.6 & 159.9 & 899.4 \\ 
		\hline 
		\textbf{example 3}& 4393.2 & 8850.6 & 49612.1 \\ 
		\hline
	\end{tabular}
	\caption{memory consumption of one macroscopic integration point in [kilobytes] (SI unit)}
	\label{memory}
\end{table}

{As a forth example, a cantilever beam with a fiber-reinforced material adopted from \citep{Tan2020} is investigated under plane stress conditions. The microstructure consists of a soft elastic-M\textsc{ises}-plastic matrix and stiff, elastic fibers. Quadratic elements with reduced integration (Abaqus \texttt{CPS8R}) are used on the macroscopic scale, while linear quadrilateral elements are employed on the microscopic level. 
\citet{Tan2020} considered this example as a benchmark problem for their \enquote{Direct-FE\textsuperscript{2} method}. This strategy is monolithic as well, but in contrast to the present method, the large coupled linear system of equations~\eqref{monolithic_linearization} for the updates of macroscopic and microscopic nodal displacements is solved in a single step (within Abaqus), without static condensation. 
The resulting load-deflection curves in Figure~\ref{fig:U_RF_beam_direct_fe2} show that both implementations yield virtually the same result, as expected. The small differences can be attributed to the differing quadrature scheme at the microscopic scale, selective reduced integration for the Direct-FE\textsuperscript{2} method compared to full integration within the present implementation. 
The last row of Figure \ref{fig:all_benchmarking_results} shows that the performance of the Direct FE\textsuperscript{2} scheme is similar to the present staggered implementation\footnote{\citet{Raju2021} even found a computational gain of the Direct FE\textsuperscript{2} method compared to a staggered implementation with finite difference quotient for the macroscopic tangent stiffness. In contrast, the present staggered implementation computes the macroscopic tangent by static condensation, thus reducing the required number of matrix factorizations at the microscopic scale.}. 
The present monolithic implementation reduces the time of computation by 47\% compared to the Direct FE\textsuperscript{2} computation if the factorized stiffness matrices are stored and by 36\% if the latter option is not used, as visualized in Figure \ref{fig:all_benchmarking_results}. The reason for the performance difference of both monolithic schemes can be found in the bandwidth $b$ of the linear systems of equations to be solved. In general, the computational cost is proportional to the number of unknowns $N$ and to $b^2$ for direct sparse matrix solvers \citep{wriggers_nonlinear_fem}. Due to the proportionality to $N$, it makes no difference with respect to the computational effort, whether a large system with $N$ unknows is to be solved (as in Direct FE\textsuperscript{2}) or if a number of small systems with the same total number of unknowns is solved (as in the present implementation with static condensation), provided the bandwidths are the same. However, in the Direct FE\textsuperscript{2}, a macroscopic node is connected to all boundary nodes of all RVEs that belong to the macroscopic elements to which the node is connected. The bandwidth with Direct FE\textsuperscript{2} is thus inevitably higher than for the decoupled systems of equations which have to be solved with the present monolithic method (or the conventional staggered scheme). This naturally has direct consequences for the respective total computational costs. The advantage of Direct FE\textsuperscript{2} lies in, of course, its simple implementation as a pure preprocessing script.}

\begin{figure}[!h]
	\centering
	\includegraphics[width=\linewidth]{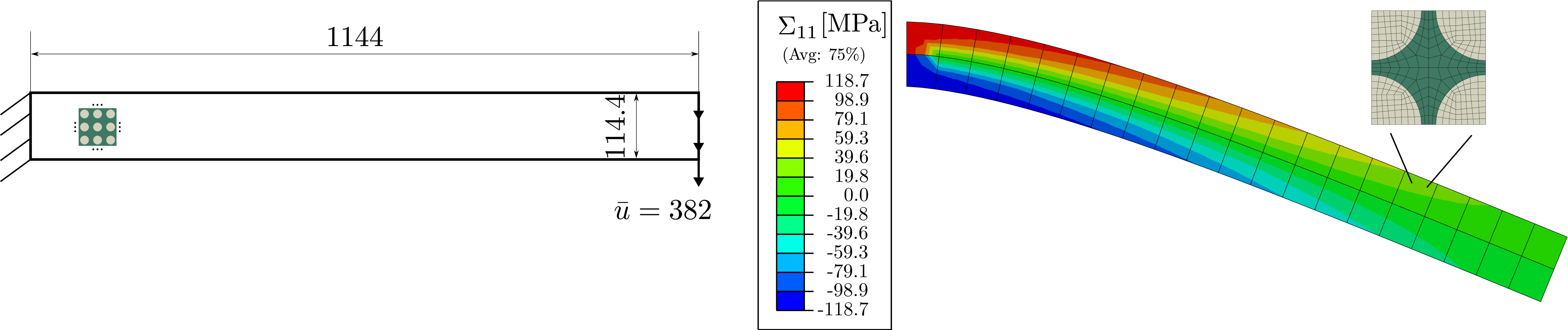}
	\caption{Example 4: Bending of a fiber-reinforced composite beam, with stiff, elastic fibres ($E=230.0$\:GPa, $\nu=0.2$) and a soft, elastic-M\textsc{ises}-plastic matrix ($E=3.9$\:GPa,\ $\nu=0.39$,\ $Y_\text{pl\,0}=29.0$\:MPa,\ $\varepsilon_\text{pl\,0}=0.0$,\ $Y_\text{pl\,1}=59.0$\:MPa,\  $\varepsilon_\text{pl\,1}=0.02$,\ $Y_\text{pl\,2}=93.0$\:MPa,\ $\varepsilon_\text{pl\,2}=0.049$,\ $Y_\text{pl\,3}=95.5$\:MPa,\ $\varepsilon_\text{pl\,0}=0.058$)}
	\label{fig:directfe2}
\end{figure}

In all benchmark cases, the proposed monolithic algorithm showed a better overall performance than the staggered solution strategy, even though a slightly higher number of macroscopic iterations was required in some time increments than with the staggered scheme. Nevertheless, the considerable reduction in the systems of equations to be solved at the microscopic scale outweighs this effect by far. Furthermore, the simulations with the staggered scheme required a few more cut-backs of time increments due to divergence at the microscale, presumably when the current macroscopic guess was outside the convergence radius of the microscopic problem.

\begin{figure}[!h]
	\centering
	\includegraphics[width=0.5\linewidth]{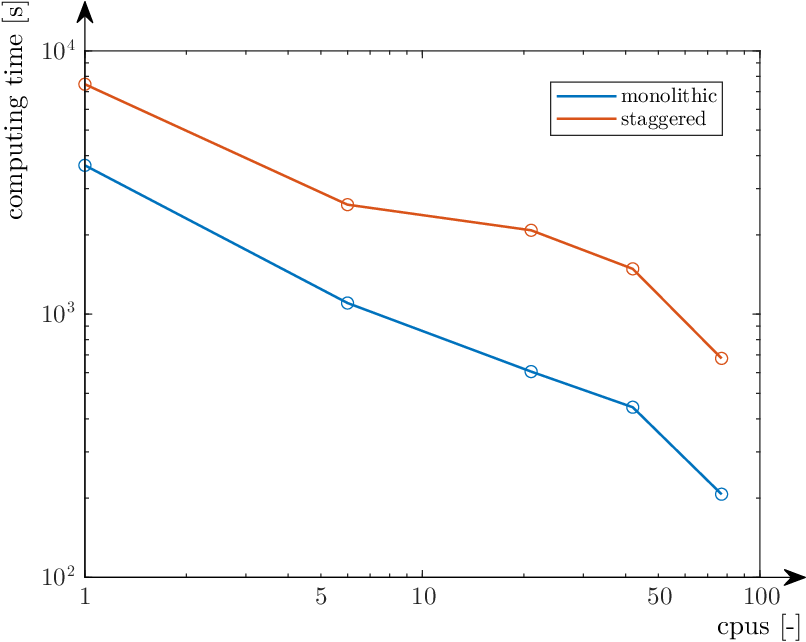}
	\caption{Computing time in dependence of number of CPU cores.}
	\label{fig:parallelization}
\end{figure}

In the previous examples, the computations were run on a single CPU core. For actual engineering applications, parallelization is of tremendous importance, as enormous amounts of computing time can be saved. 
In general, FE\textsuperscript{2} methods are very well suited for parallelization, since the microscale problems can be solved independently of each other, see e.g.\ \cite{Geers_multiscale_2017,Nakshatrala_multiscale_topology_2013}).
%Since the present implementation is fully parallelizable, it is possible to compute multiple micro-scale FE problems simultaneously. 
In order to investigate the scalability, benchmark example 3 (using the small deformation theory) was run with 1, 6, 21, 42 and 77 CPU cores over multiple compute nodes (distributed memory architecture), using the monolithic algorithm with stored factorization and the staggered scheme, respectively. The results in Figure~\ref{fig:parallelization} demonstrate a very good overall scalability. 
Generally, the monolithic algorithm has a higher potential of profiting from parallelization, since the effort for each macroscopic integration point, namely assembling and solving one system of equations for the update of the microscale displacements, is approximately the same. In contrast, the number of necessary microscopic iterations and thus the required computational time may differ considerably between the macroscopic integration points in the staggered scheme, so that it is more likely that one CPU core has to wait for the others to finish their microscale iterations.

%--------------------------------------------
\section{Conclusions}
%--------------------------------------------
The computational costs of FE\textsuperscript{2} simulations are generally very high. Currently this fact, and the lack of implementations in commercial FE software, often hinder their application in actual engineering problems. Hence, a monolithic solution strategy to FE\textsuperscript{2} problems has been outlined in the present contribution, which solves the macro- and the micro-scale problems in a common Newton-Raphson iteration loop. It was shown that the microscopic degrees of freedom can be removed from the global system of equations through static condensation, due to their weak coupling. The main idea behind the monolithic approach is to lower the number of costly microscopic iterations considerably and therefore the total computational effort.

The presented theory was implemented into Abaqus through the \texttt{UMAT} interface in a concise, flexible and fully parallelizable manner. The program has been successfully verified through many comparative calculations.

The differences in computational effort between the staggered and monolithic algorithms depend on many factors, of which the most significant is the nonlinear character of the considered problems, as the monolithic and staggered schemes are equivalent in linear cases. To quantify the differences, several benchmark tests were carried out, which demonstrated that the monolithic approach is capable of saving vast amounts of computing time. More specifically, the benchmark examples showed that up to 60\% of the computational effort were saved. To also examine the parallelizability of the FE\textsuperscript{2}-method in practice, a parallelization study was performed and its results indicated a good overall scalability. The monolithic Direct-FE\textsuperscript{2} method was compared with the proposed approach. The proposed monolithic algorithm with static condensation turned out to be faster, due to a lower bandwidth of the resulting systems of equations. The Direct-FE\textsuperscript{2}, on the other hand, is method is more flexible due to the availability of the Abaqus features on both scales.

It is therefore concluded that the presented monolithic approach --- particularly when combined with massive parallelization --- can reduce the computing time to a level at which it becomes realistic to employ the FE\textsuperscript{2}-method in simulations of real-life engineering applications, thereby enhancing such analyses by multi-scale information.

\section*{Acknowledgment}
The authors thank Vincent Tan and Karthikayen Raju for providing the Direct FE\textsuperscript{2} implementation (Example 4) and for many fruitful discussions.
Furthermore, the funding by the Deutsche Forschungsgemeinschaft (DFG, German Research Foundation) via the SFB 920 \enquote{Multi-Functional Filters for Metal Melt Filtration --- A Contribution towards Zero Defect Materials} 
-- project ID 169148856 -- is gratefully acknowledged.

% Nur die bib-Vorlage von Elsevier (bei vielen Latex-distributionen schon dabei)
\bibliographystyle{elsarticle-harv}
\bibliography{literatur_FEmonolith}

\end{document}